\documentclass[12pt]{article}
\usepackage{amssymb,amsfonts,amsmath, psfrag,eepic,colordvi}
\newcommand{\rmnum}[1]{\romannumeral #1}
\topskip  
\numberwithin{equation}{section}
\newtheorem{theorem}{Theorem}[section]

\newtheorem{remark}[theorem]{Remark}
\newtheorem{lemma}[theorem]{Lemma}

\newtheorem{example}[theorem]{Example}

\begin{document}
\parskip 6pt

\pagenumbering{arabic}
\def\sof{\hfill\rule{2mm}{2mm}}
\def\ls{\leq}
\def\gs{\geq}
\def\SS{\mathcal S}
\def\qq{{\bold q}}
\def\MM{\mathcal M}
\def\TT{\mathcal T}
\def\EE{\mathcal E}
\def\lsp{\mbox{lsp}}
\def\rsp{\mbox{rsp}}
\def\pf{\noindent {\it Proof.} }
\def\mp{\mbox{pyramid}}
\def\mb{\mbox{block}}
\def\mc{\mbox{cross}}
\def\qed{\hfill \rule{4pt}{7pt}}
\def\block{\hfill \rule{5pt}{5pt}}

\begin{center}
{\Large\bf  Bijective counting of humps and peaks in  $(k,a)$-paths  }
\vskip 6mm
\end{center}

\begin{center}
{\small   Sherry H. F. Yan \\[2mm]
 Department of Mathematics, Zhejiang Normal University, Jinhua
321004, P.R. China
\\[2mm]
 huifangyan@hotmail.com
  \\[0pt]
}
\end{center}

\noindent {\bf Abstract.} Recently, Mansour and Shattuck   related
the total number of humps  in all of the $(k, a)$-paths of order $n$
to the number of super $(k, a)$-paths, which generalized   previous
results concerning the cases when $k = 1$ and $a = 1$ or $a =
\infty$.  They also derived  a relation on the total number of peaks
in all of the $(k, a)$-paths of order $n$ and the number of super
$(k, a)$-paths, and asked for bijective proofs. In this paper, we
will give bijective proofs of these two relations.




\section{Introduction}
 A $(k,a)$-path of order $n$ is a lattice path in $\emph{Z}\times \emph{Z}$ from $(0,0)$ to $(n,0)$
 using up steps $(1,k)$, down steps $(1,-1) $and horizontal steps $(a,0)$ and never
  lying below the $x$-axis.
   Denote by $\mathcal{P}_n(k,a)$ the set of all $(k,a)$-paths of order $n$.
  Note that $\mathcal{P}_n(k, \infty)$, $\mathcal{P}_n(1,\infty)$, $\mathcal{P}_n(1,1)$, and $P_n(1,2)$ are the
   set of $k$-ary paths \cite{DM}, Dyck paths, Motzkin paths and Schr\"{o}der paths, respectively.
   If a $(k,a)$-path of order $n$ is allowed to go below the $x$-axis, then it is called a {\em super} $(k,a)$-path of order $n$. Denote by $\mathcal{SP}_n(k,a)$ the set of all super $(k,a)$-paths of order $n$.   A {\em peak} in a $(k,a)$-path is an up step followed by a down step. A {\em hump} in a $(k,a)$-path is an up step followed by zero or more horizontal steps followed by a down step. We denote by $\# Peaks(P)$ and $\# Humps(P)$ the number of peaks and humps in a path $P$.

 Using a recurrence relation and the WZ method, Regev \cite{R} proved that

 \begin{equation}\label{hump1}
 2\sum_{P\in \mathcal{P}_n(1,1)}\# Humps(P)=|\mathcal{SP}_{n}(1,1)|-1.
  \end{equation}

 \begin{equation}\label{peak1}
 2\sum_{P\in \mathcal{P}_n(1,\infty)}\#Peaks(P)=|\mathcal{SP}_{n}(1,\infty)|
  \end{equation}

 The bijective proofs of Formulae (\ref{hump1}) and (\ref{peak1}) were given by Ding and Du \cite{DD}. They also derived the following analogous result for Schr\"{o}der  paths
  \begin{equation}\label{schroder}
 2\sum_{P\in \mathcal{P}_n(1,2)}\#Humps(P)=|\mathcal{SP}_{n}(1,2)|-1
  \end{equation}

Recently, Mansour and Shattuck  \cite{MS}  proved that
\begin{equation}\label{hump}
 (k+1)\sum_{P\in \mathcal{P}_n(k,a)}\# Humps(P)=|\mathcal{SP}_{n}(k,a)|-\delta_{a|n},
  \end{equation}

 \begin{equation}\label{peak}
 (k+1)\sum_{P\in \mathcal{P}_n(k,a)}\# Peaks(P)=|\mathcal{SP}_{n}(k,a)|-|\mathcal{SP}_{n-a}(k,a)|,
  \end{equation}
 where $\delta_{a|n}=1$ if $a$ divides $n$ or $0$ otherwise, and asked for bijective proofs.    Specializing $k=1$ and $a=1$, $a=\infty$ or $a=2$ in Formulae (\ref{hump}) and (\ref{peak}) gives
 Formulae (\ref{hump1}), (\ref{peak1}) and (\ref{schroder}). The main objective of this paper  is to give bijective proofs of Formulae (\ref{hump}) and (\ref{peak}) in answer to the problem posed by Mansour and Shattuck.
 As a consequence, our bijection also allows us to get the
 enumeration of $k$-ary paths with respect to the number of peaks.

 \section{The bijective proofs}
 In this section, we will give bijective proofs of (\ref{hump}) and (\ref{peak}). We begin with some necessary definitions and notations.

    Throughout this paper we identify a   path with a word  by encoding each up step by the letter $U$, each down step by $D$ and each horizontal step by $H$.    Let $p$ be a step  running from the point $(x_1, y_1)$ to the point $(x_2, y_2)$. Then we say that the point $(x_1,y_1)$ is its  {\em starting
   } point, and the  point $(x_2, y_2)$ is its {\em ending} point. The {\em starting} point and the {\em ending } point  of a path are defined analogously.
       A point $(x,y)$
    of a lattice path $P$
    is said to be a {\em return } point   if $y=0$ and $(x,y)$ is not the starting point of
    $P$.    An up step is said to {\em intersect} the $x$-axis if its starting point lies weakly below the $x$-axis and its ending point lies weakly above the $x$-axis.
    If $P=p_1p_2\ldots p_n$ is a path,  then  the {\em reverse} of the path,
denoted by $\widehat{P}$, is defined  by $p_np_{n-1}\ldots p_{1}$.
For example, the reverse of the path  $P=HUDDUDH$ is given by  $HDUDDUH$.

  A $(k,a)$-path is said to be   {\em $m$-peak} (resp. {\em $m$-hump})  colored if  exactly one   peak  (resp. hump )
  is assigned  by any of the $m$ colors.  Denote by $\mathcal{PP}^m_n(k,a)$  (resp.  $\mathcal{PM}^m_n(k,a)$ ) the set of all $m$-peak (resp. $m$-hump) colored  $(k,a)$-paths.

 Observe that the left-hand sides of (\ref{hump}) and $(\ref{peak})$ are equal to $|\mathcal{PM}^{k+1}_n(k,a)|$
  and $|\mathcal{PP}^{k+1}_n(k,a)|$, respectively. Moreover,
 the right-hand side of  (\ref{hump}) counts the number of super $(k,a)$-paths of order $n$ and  with at least one up step.  The right-hand side of  (\ref{peak}) counts the number of super $(k,a)$-paths of order $n$  which do not
  start  with horizontal steps. Denote by $\mathcal{S'}_n(k,a)$ and $\mathcal{S''}_n(k,a)$ the set
  of all super $(k,a)$-paths of order $n$ and  with at least one up step, and the set of all super $(k,a)$-paths of order $n$
  which do not start with horizontal steps, respectively. Thus Formulae  (\ref{hump}) and (\ref{peak}) can be rewritten as
 \begin{equation}\label{ref1}
 | \mathcal{PM}^{k+1}_n(k,a)|=|\mathcal{S'}_n(k,a)|,
 \end{equation}

  \begin{equation}\label{ref2}
 | \mathcal{PP}^{k+1}_n(k,a)|=|\mathcal{S''}_n(k,a)|.
 \end{equation}
 In order to prove Formulae (\ref{ref1} ) and (\ref{ref2} ), we will establish   a bijection  $\phi$  between the set $\mathcal{PM}^{k+1}_n(k,a)$ and the set $\mathcal{S'}_n(k,a)$.  Moreover,   we show that the map $\phi$   restricted to the set $\mathcal{PP}^{k+1}_n(k,a)$ gives a bijection between   the set $\mathcal{PP}^{k+1}_n(k,a)$ and the set $\mathcal{S''}_n(k,a)$.

 Now we  proceed to describe the map $\phi$ from the set $\mathcal{PM}^{k+1}_n(k,a)$ to the set
 $\mathcal{S}'_n(k,a)$. Let $P=p_1p_2\ldots p_n$ be   a $(k+1)$-hump colored
 $(k,a)$-path of order $n$.
 Suppose that the hump  $p_lp_{l+1}\ldots p_{m}=UH^{m-l-1}D$ is colored by $c$,  where $H^{i}$ denotes $i$ consecutive horizontal steps.
Assume that $p_l$ goes from the point  $(x_1, h)$ to the point
$(x_1+1, h+k)$ for some nonnegative integers $x_1$ and    $h$.

 Then  the  path  $P$ can be uniquely decomposed as
 $$
  R_1 P' p_l H^{m-l-1} d_1 R_2 d_2\ldots R_{k}d_k  P'',
$$
where
\begin{itemize}
\item each $d_i$ is the first  step after  $p_l$ that goes from the line $y=h+k+1-i$ to the line $y=h+k-i$;
 \item $P'$ is the (possibly empty) section of $P$       which is to the left of $p_l$,  starts with an up step,  and lies   strictly above the $x$-axis except for the starting point;
     \item  each $R_i$ is a (possibly empty) $(k,a)$-path;
    \item $P''$ is the remaining  section  of $P$ after $d_k$.
  \end{itemize}
  Obviously, each $d_i$ is a down step and the subpath $P''$ goes from the line $y=h$ to the $x$-axis in $P$.
Now we proceed to construct $\phi(P)$ as follows:
\begin{itemize}
\item[$(\rmnum{1})$]  if $c=1$, then set
$$
\phi(P)=H^{m-l-1} p_l R_{1}d_{1}R_2d_2 \ldots R_k d_k P'' P';
$$
\item[$(\rmnum{2})$]  if $c=k+1$, then set
$$
\phi(P)=H^{m-l-1} d_1  \widehat{R_1}  d_2 \widehat{R_2}  \ldots
d_k \widehat{R_k } p_l P'' P';
$$
\item[$(\rmnum{3} )$] if $1<c<k+1$, then set
$$\phi(P)=H^{m-l-1} d_1  \widehat{R_1}  d_2  \widehat{R_2}  \ldots d_{c-1}  \widehat{R_{c-1}}  p_l R_{c}d_{c} \ldots R_k d_k P'' P'.$$
\end{itemize}

  According to the construction of the map $\phi$,
  we preserve the number of up steps, the number of down steps and the number of horizontal steps.
   Moreover, there is at least one up step in the resulting path. Hence,  the map $\phi$ is well defined, that is,  $\phi(P)\in \mathcal{S'}_n(k,a)$.

\begin{remark}
Our map $\phi$ restricted to case when $k=1$  is different from the bijection given by Ding and Du \cite{DD}.
\end{remark}

\begin{figure}[h,t]
\begin{center}
\begin{picture}(120,20)
\setlength{\unitlength}{3mm} \linethickness{0.4pt}
\put(0,0){\line(1,0){38}} \put(0,0){\circle*{0.2}}
 \put(0,0){\line(1,3){1}}
 \put(1,3){\circle*{0.2}}
 \put(1,3){\line(1,-1){1}}\put(2,2){\line(1,-1){1}}\put(2,2){\circle*{0.2}}\put(3,1){\circle*{0.2}}
 \put(3,1){\line(1,-1){1}}\put(4,0){\circle*{0.2}}
 \put(4,0){\circle*{0.2}}
 \put(4,0){\line(1,3){1}}
 \put(5,3){\circle*{0.2}}
 \put(5,3){\line(1,-1){1}} \put(6,2){\circle*{0.2}}
 \put(6,2){\line(1,3){1}} \put(7,5){\circle*{0.2}}
 \put(7,5){\line(1,-1){1}} \put(8,4){\circle*{0.2}}
 \put(8,4){\line(1,-1){1}} \put(9,3){\circle*{0.2}}
 \put(9,3){\line(1,-1){1}} \put(10,2){\circle*{0.2}}
 \put(10,2){\line(1,3){1}} \put(11,5){\circle*{0.2}}
 \put(11,5){\line(1,-1){1}} \put(12,4){\circle*{0.2}}
  \put(12,4){\line(1,3){1}} \put(13,7){\circle*{0.2}}
   \put(13,7){\line(2,0){2}} \put(15,7){\circle*{0.2}}
  \put(13,7){*}
  \put(15,7){\line(1,-1){1}} \put(16,6){\circle*{0.2}}
  \put(16,6){\line(1,3){1}} \put(17,9){\circle*{0.2}}
  \put(17,9){\line(1,-1){1}} \put(18,8){\circle*{0.2}}
  \put(18,8){\line(1,-1){1}} \put(19,7){\circle*{0.2}}
  \put(19,7){\line(1,-1){1}} \put(20,6){\circle*{0.2}}

  \put(20,6){\line(1,-1){1}} \put(21,5){\circle*{0.2}}
  \put(21,5){\line(1,-1){1}} \put(22,4){\circle*{0.2}}
  \put(22,4){\line(1,-1){1}} \put(23,3){\circle*{0.2}}
  \put(23,3){\line(1,3){1}} \put(24,6){\circle*{0.2}}
  \put(24,6){\line(1,-1){1}} \put(25,5){\circle*{0.2}}
  \put(25,5){\line(1,-1){1}} \put(26,4){\circle*{0.2}}
  \put(26,4){\line(1,-1){1}} \put(27,3){\circle*{0.2}}
  \put(27,3){\line(1,-1){1}} \put(28,2){\circle*{0.2}}

    \put(28,2){\line(1,3){1}} \put(29,5){\circle*{0.2}}
  \put(29,5){\line(1,-1){1}} \put(30,4){\circle*{0.2}}
  \put(30,4){\line(1,-1){1}} \put(31,3){\circle*{0.2}}
  \put(31,3){\line(1,-1){1}} \put(32,2){\circle*{0.2}}
  \put(32,2){\line(1,-1){1}} \put(33,1){\circle*{0.2}}
   \put(33,1){\line(1,-1){1}} \put(34,0){\circle*{0.2}}
   \put(0,9){\line(1,0){38}}
   \put(0,0){\line(0,1){9}}\put(4,0){\line(0,1){9}}
   \put(1.5,9.5){$R_1$}
\put(21,0){\line(0,1){9}} \put(21.2,9.5){\small
$d_3$}\put(19.8,9.5){\small$d_2$}
   \put(7,9.5){$P'$}
    \put(12,0){\line(0,1){9}}
\put(15,0){\line(0,1){9}}
    \put(16,0){\line(0,1){9}}\put(20,0){\line(0,1){9}}
   \put(17.5,9.5){$R_2$}
   \put(22,0){\line(0,1){9}}
   \put(30.5,9.5){$P''$}
\put(13,0){\line(0,1){9}} \put(11.8, 9.5){$p_l$}
   \put(34,0){\line(1,3){1}}
    \put(35,3){\circle*{0.2}}
    \put(35,3){\line(1,-1){1}}
    \put(36,2){\circle*{0.2}}
    \put(36,2){\line(1,-1){1}}
    \put(37,1){\circle*{0.2}}
   \put(37,1){\line(1,-1){1}}
    \put(38,0){\circle*{0.2}}
    \put(38,0){\line(0,1){9}}
\put(14.8, 9.5){\small$d_1$}
\end{picture}
\end{center}
\caption{The decomposition of a $(3,2)$-path.} \label{fig1}
\end{figure}
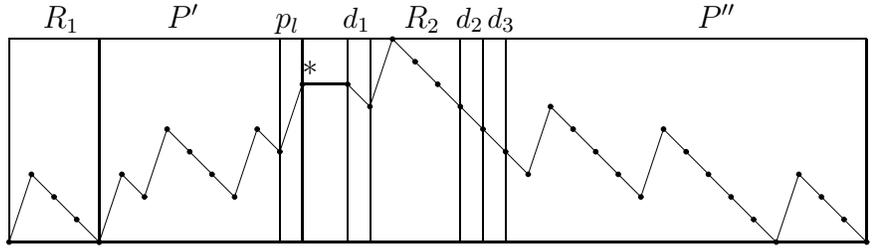

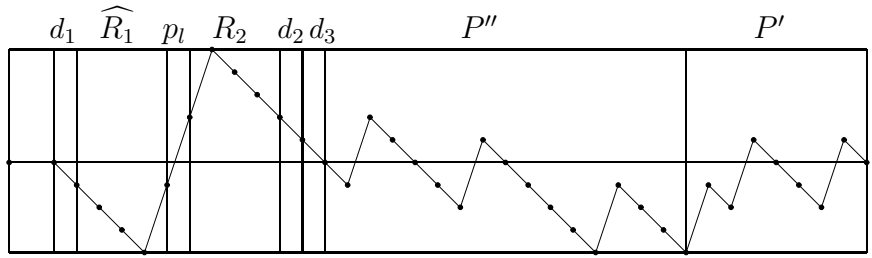
\begin{figure}[h,t]
\begin{center}
\begin{picture}(120,30)
\setlength{\unitlength}{3mm} \linethickness{0.4pt}
\put(0,0){\line(0,1){9}}\put(38,0){\line(0,1){9}}
\put(0,4){\line(1,0){38}}\put(0,9){\line(1,0){38}}\put(0,0){\line(1,0){38}}
  \put(0,
4){\line(1,0){2}} \put(0,4){\circle*{0.2}}\put(2,4){\circle*{0.2}}
\put(2, 4){\line(1,-1){1}} \put(3,3){\circle*{0.2}} \put(3,
3){\line(1,-1){1}} \put(4,2){\circle*{0.2}} \put(4,
2){\line(1,-1){1}} \put(5,1){\circle*{0.2}} \put(5,
1){\line(1,-1){1}} \put(6,0){\circle*{0.2}} \put(6,
0){\line(1,3){1}} \put(7,3){\circle*{0.2}}\put(7, 3){\line(1,3){1}}
\put(8,6){\circle*{0.2}}\put(8, 6){\line(1,3){1}}
\put(9,9){\circle*{0.2}}

\put(9, 9){\line(1,-1){1}} \put(10,8){\circle*{0.2}} \put(10,
8){\line(1,-1){1}} \put(11,7){\circle*{0.2}} \put(11,
7){\line(1,-1){1}} \put(12,6){\circle*{0.2}}

\put(12, 6){\line(1,-1){1}} \put(13,5){\circle*{0.2}} \put(13,
5){\line(1,-1){1}} \put(14,4){\circle*{0.2}} \put(14,
4){\line(1,-1){1}} \put(15,3){\circle*{0.2}}\put(15,
3){\line(1,3){1}} \put(16,6){\circle*{0.2}}

\put(16, 6){\line(1,-1){1}} \put(17,5){\circle*{0.2}} \put(17,
5){\line(1,-1){1}} \put(18,4){\circle*{0.2}} \put(18,
4){\line(1,-1){1}} \put(19,3){\circle*{0.2}}\put(19,
3){\line(1,-1){1}} \put(20,2){\circle*{0.2}}

\put(20, 2){\line(1,3){1}} \put(21,5){\circle*{0.2}}

\put(21, 5){\line(1,-1){1}} \put(22,4){\circle*{0.2}} \put(22,
4){\line(1,-1){1}} \put(23,3){\circle*{0.2}} \put(23,
3){\line(1,-1){1}} \put(24,2){\circle*{0.2}}\put(24,
2){\line(1,-1){1}} \put(25,1){\circle*{0.2}}\put(25,
1){\line(1,-1){1}} \put(26,0){\circle*{0.2}}

\put(26, 0){\line(1,3){1}} \put(27,3){\circle*{0.2}}

\put(27, 3){\line(1,-1){1}} \put(28,2){\circle*{0.2}} \put(28,
2){\line(1,-1){1}} \put(29,1){\circle*{0.2}} \put(29,
1){\line(1,-1){1}} \put(30,0){\circle*{0.2}}

\put(30, 0){\line(1,3){1}} \put(31,3){\circle*{0.2}}

\put(31, 3){\line(1,-1){1}} \put(32,2){\circle*{0.2}}

\put(32, 2){\line(1,3){1}} \put(33,5){\circle*{0.2}}

\put(33, 5){\line(1,-1){1}} \put(34,4){\circle*{0.2}} \put(34,
4){\line(1,-1){1}} \put(35,3){\circle*{0.2}} \put(35,
3){\line(1,-1){1}} \put(36,2){\circle*{0.2}}

\put(36, 2){\line(1,3){1}} \put(37,5){\circle*{0.2}}
\put(2,0){\line(0,1){9}} \put(1.8, 9.5){\small$d_1$}
\put(13,0){\line(0,1){9}} \put(11.9, 9.5){\small$d_2$} \put(13.3,
9.5){\small$d_3$}
 \put(37, 5){\line(1,-1){1}}
\put(38,4){\circle*{0.2}}
\put(30,0){\line(0,1){9}}\put(14,0){\line(0,1){9}}
\put(3,0){\line(0,1){9}}\put(7,0){\line(0,1){9}}\put(8,0){\line(0,1){9}}\put(12,0){\line(0,1){9}}
\put(4, 9.5){$\widehat{ R_1 }$}\put(9, 9.5){$  R_2  $}\put(33,
9.5){$ P'  $}\put(20, 9.5){$ P''  $}\put(6.8, 9.5){$ p_l  $}
\end{picture}
\end{center}
\caption{The application of $\phi$ to the path shown in Figure \ref{fig1}.} \label{fig2}
\end{figure}

 \begin{example}

An example of  the decomposition of a $(3,2)$-path is shown in Figure \ref{fig1}, where
  the colored hump is marked by a star and $R_3$ is an empty path.
    Suppose that the hump is colored by $2$. By applying the map $\phi$, we get its corresponding super $(3,2)$-path   shown in Figure \ref{fig2}.
\end{example}

In order to show that $\phi$ is a bijection, we  describe a map
$\psi$ from the set $\mathcal{S'}_n(k,a)$  to the set
$\mathcal{PM}^{k+1}_n(k,a)$. Given a super  $(k,a)$-path  $Q=q_1q_2\ldots q_n\in
\mathcal{S'}_n(k,a)$,  let  $q_l$ be  the leftmost up step that
intersects the $x$-axis. Suppose that $q_l$ goes from the point
$(x_1, p)$ to the point $(x_1+1, q)$, where $q-p=k$. Let $A$  be  the  first return point to the right of the point  $(x_1, p)$ and $B$ be the
the  lowest point to the right of    $(x_1,p)$ in $Q$. If there are more than one such lowest point,  we choose $B$ to be  the rightmost one. Then we generate $\psi(Q)$ as follows.
 \begin{itemize}

   \item[$(\rmnum{1}') $]  If $q>0$ and $p=0$, then $Q$ can be uniquely decomposed as
   $$
  Q=H^{l-1} q_l R_1 d_1 R_2 d_2\ldots  R_k d_k  Q' Q'',
  $$
where
  \begin{itemize}

  \item  each $d_i$  is the first step that goes from the line $y=k+1-i$ to the line $y=k-i$;
  \item  each $R_i$ is a (possibly empty) $(k,a)$-path;
  \item   $Q'$  is the section
of $Q$ which goes from the point $A$ to the point $B$;
  \item  $Q''$ is the remaining section of $Q$.
  \end{itemize}

   Then
  set
  $$\psi(Q)=R_1 Q'' q_l H^{l-1} d_1 R_2 d_2\ldots R_kd_k Q',$$
  where the hump $q_l H^{l-1} d_1$ is colored by $1$.

\item[$(\rmnum{2}')$ ] If $q=0$, then $Q$ can be uniquely decomposed as
  $$
  Q= H^{m}    d_1  \widehat{R_1}  d_2 \widehat{ R_2}  \ldots    d_k  \widehat{R_k}  q_l Q'Q'',
  $$
   where
\begin{itemize}
\item $m$ is an nonnegative integer;
\item each $d_i$ is the last step to the left of the point $(x_1, p)$   that  goes form the line $y=-i+1$ to the line $y=-i$;
\item each $R_i$ is a (possibly empty) $(k,a)$-path;
 \item   $Q'$  is the section
of $Q$ which goes from the point $A$ to the point $B$;
\item $Q''$ is the remaining section of $Q$.
\end{itemize}

  Then set
  $$\psi(Q)=R_1Q'' q_l H^{m} d_1 R_2 d_2\ldots R_kd_k Q',$$
  where the hump $q_l H^{m} d_1$ is colored by $k+1$.

  \item[$(\rmnum{3}') $]  If $q>0$ and $p<0$, then $Q$ can be uniquely decomposed
  as
   $$
  Q= H^{m}    d_1   \widehat{R_1}   d_2  \widehat{ R_2}   \ldots   d_{|p|}    \widehat{R_{|p|} }  q_l R_{|p|+1}d_{|p|+1} \ldots R_k d_k Q' Q'',
  $$
where
\begin{itemize}
\item $m$ is an nonnegative integer;
\item for $1\leq i\leq |p|$, each $d_i$ is the last step to the  left of  $q_l$ that goes from the line $y=-i+1$ to the line $y=-i$;
\item for $|p|+1\leq i\leq k $, each $d_i$ is the first step   that goes from the line $y=k+1-i$ to the line $y=k-i$;
\item   each $R_i$ is a (possibly empty) $(k,a)$-path;
 \item   $Q'$  is the section
of $Q$ which goes from the point $A$ to the point $B$;
\item $Q''$ is the remaining section of $Q$.
\end{itemize}
 Then set
  $$\psi(Q)=  R_1  Q'' q_l H^{m} d_1  R_2  d_2\ldots  R_kd_k Q',$$
  where the hump $q_l H^{m} d_1$ is colored by $|p|+1$.
 \end{itemize}

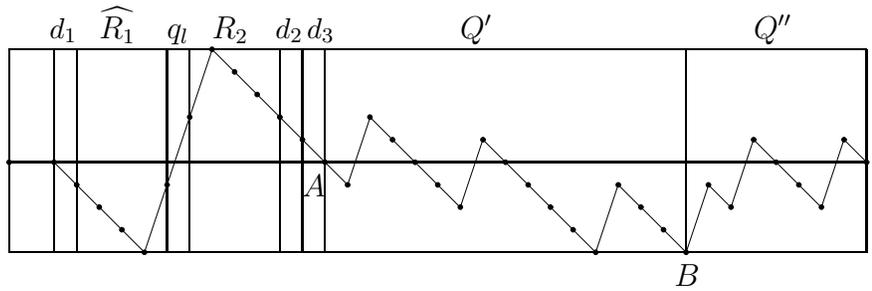
\begin{figure}[h,t]
\begin{center}
\begin{picture}(120,30)
\setlength{\unitlength}{3mm} \linethickness{0.4pt}
\put(0,0){\line(0,1){9}}\put(38,0){\line(0,1){9}}
\put(0,4){\line(1,0){38}}\put(0,9){\line(1,0){38}}\put(0,0){\line(1,0){38}}
  \put(0,
4){\line(1,0){2}} \put(0,4){\circle*{0.2}}\put(2,4){\circle*{0.2}}
\put(2, 4){\line(1,-1){1}} \put(3,3){\circle*{0.2}} \put(3,
3){\line(1,-1){1}} \put(4,2){\circle*{0.2}} \put(4,
2){\line(1,-1){1}} \put(5,1){\circle*{0.2}} \put(5,
1){\line(1,-1){1}} \put(6,0){\circle*{0.2}} \put(6,
0){\line(1,3){1}} \put(7,3){\circle*{0.2}}\put(7, 3){\line(1,3){1}}
\put(8,6){\circle*{0.2}}\put(8, 6){\line(1,3){1}}
\put(9,9){\circle*{0.2}}

\put(9, 9){\line(1,-1){1}} \put(10,8){\circle*{0.2}} \put(10,
8){\line(1,-1){1}} \put(11,7){\circle*{0.2}} \put(11,
7){\line(1,-1){1}} \put(12,6){\circle*{0.2}}

\put(12, 6){\line(1,-1){1}} \put(13,5){\circle*{0.2}} \put(13,
5){\line(1,-1){1}} \put(14,4){\circle*{0.2}} \put(14,
4){\line(1,-1){1}} \put(15,3){\circle*{0.2}}\put(15,
3){\line(1,3){1}} \put(16,6){\circle*{0.2}}

\put(16, 6){\line(1,-1){1}} \put(17,5){\circle*{0.2}} \put(17,
5){\line(1,-1){1}} \put(18,4){\circle*{0.2}} \put(18,
4){\line(1,-1){1}} \put(19,3){\circle*{0.2}}\put(19,
3){\line(1,-1){1}} \put(20,2){\circle*{0.2}}

\put(20, 2){\line(1,3){1}} \put(21,5){\circle*{0.2}}

\put(21, 5){\line(1,-1){1}} \put(22,4){\circle*{0.2}} \put(22,
4){\line(1,-1){1}} \put(23,3){\circle*{0.2}} \put(23,
3){\line(1,-1){1}} \put(24,2){\circle*{0.2}}\put(24,
2){\line(1,-1){1}} \put(25,1){\circle*{0.2}}\put(25,
1){\line(1,-1){1}} \put(26,0){\circle*{0.2}}

\put(26, 0){\line(1,3){1}} \put(27,3){\circle*{0.2}}

\put(27, 3){\line(1,-1){1}} \put(28,2){\circle*{0.2}} \put(28,
2){\line(1,-1){1}} \put(29,1){\circle*{0.2}} \put(29,
1){\line(1,-1){1}} \put(30,0){\circle*{0.2}}

\put(30, 0){\line(1,3){1}} \put(31,3){\circle*{0.2}}

\put(31, 3){\line(1,-1){1}} \put(32,2){\circle*{0.2}}

\put(32, 2){\line(1,3){1}} \put(33,5){\circle*{0.2}}

\put(33, 5){\line(1,-1){1}} \put(34,4){\circle*{0.2}} \put(34,
4){\line(1,-1){1}} \put(35,3){\circle*{0.2}} \put(35,
3){\line(1,-1){1}} \put(36,2){\circle*{0.2}}

\put(36, 2){\line(1,3){1}} \put(37,5){\circle*{0.2}}

\put(37, 5){\line(1,-1){1}} \put(38,4){\circle*{0.2}}
\put(30,0){\line(0,1){9}}\put(14,0){\line(0,1){9}}
\put(3,0){\line(0,1){9}}\put(7,0){\line(0,1){9}}\put(8,0){\line(0,1){9}}\put(12,0){\line(0,1){9}}
\put(4, 9.5){$\widehat{R_1}$}\put(9, 9.5){$  R_2  $}\put(33, 9.5){$ Q''
$}\put(20, 9.5){$ Q'  $} \put(29.5, -1.5){$B$}\put(13, 2.5){$A$}
\put(2,0){\line(0,1){9}} \put(1.8, 9.5){\small$d_1$}
\put(13,0){\line(0,1){9}} \put(11.8, 9.5){\small$d_2$} \put(13.2,
9.5){\small$d_3$}\put(7, 9.5){$q_l$}
\end{picture}
\end{center}
\caption{The decomposition of  a  super $(3,2)$-path $Q$. } \label{fig3}
\end{figure}

\begin{example}
The decomposition of a  super $(3,2)$-path $Q$ is illustrated in Figure
\ref{fig3}, where $A=(14,0)$ and $B=(30,-4)$,  the path $R_3$ is empty,  and the leftmost up step
that intersects the $x$-axis goes from the point $(7,-1)$ to the
point $(8, 2)$.  By applying the map $\psi$, we get a $(3,2)$-path $\psi(Q)$ shown  in Figure \ref{fig1}.
\end{example}

Obviously, each $d_i$ is a down step. It is easy to check that the map $\psi$ is well defined, that is,
$\psi(Q)\in \mathcal{PM}^{k+1}_n(k,a) $. From the construction of
the map $\phi$, we see that  $p_l$ is the leftmost up step that
intersects the $x$-axis and the ending point of $P''$ is the last
lowest point   to the  right of  the starting point of $p_l$ in $\phi(p)$. Moreover, the starting point of $P''$ is the first return point to the right of  the starting point of $p_l$ in $\phi(P)$.
Thus, one can easily  verify that $(\rmnum{1}') $, $(\rmnum{2}') $
and $(\rmnum{3}') $ respectively reverse the procedures of
$(\rmnum{1} ) $, $(\rmnum{2} ) $ and $(\rmnum{3} ) $. This implies
that  the maps $\phi$ and $\psi$ are inverses of each other. Hence,
the map $\phi$ is a bijection.

 \begin{theorem}\label{th1}
 The map $\phi$ is a bijection between the set $\mathcal{PM}^{k+1}_n(k,a)$ and  the set $\mathcal{S}'_n(k,a)$.
 Moreover, for any $P\in \mathcal{PM}^{k+1}_n(k,a)$ whose colored hump  consists of $UH^lD$,   the corresponding
 super $(k,a)$-path  $\phi(P)$ starts with exactly $l$ consecutive horizontal steps.
 \end{theorem}

 From Theorem \ref{th1},  it follows that the bijection $\phi$ restricted to the set $\mathcal{PP}^{k+1}_n(k,a)$
  reduces to a bijection between the set $\mathcal{PP}^{k+1}_n(k,a)$ and  the set $\mathcal{S}''_n(k,a)$.
   Thus we obtain  bijective proofs of Formulae (\ref{hump}) and (\ref{peak}).

Our bijection $\phi$ also allows us to enumerate $k$-ary paths with
respect to the number of peaks. Let $\mathcal{Q}_n(k,m)$ be the set
of $k$-ary paths with $n$ up steps and  $m$ peaks in which exactly
one peak is colored by $1$. Denote by $\mathcal{S}^{UU}_n(k, m)$ the
set of super $k$-ary paths   with $n$ up steps and  $m$ peaks which
start  with at least two consecutive up steps. Denote by
$\mathcal{S}^{UD}_n(k, m)$ the set of super $k$-ary paths with $n$
up steps and  $m$ peaks which start with an up step followed
immediately by a down step.

 From the construction of the bijection $\phi$, it is easily seen that
  the bijection $\phi$ restricted to the set
$\mathcal{Q}_n(k,m)$ reduces to a bijection between the set
$\mathcal{Q}_n(k,m)$ and the set $\mathcal{S}^{UU}_n(k, m-1)\cup
\mathcal{S}^{UD}_n(k, m)$. In order to get the enumeration of
$k$-ary paths with respect to the number of peaks, we need the
following lemma.

\begin{lemma}
For $n, m\geq 1$, we have
\begin{equation}\label{c1}
|\mathcal{S}^{UU}_n(k, m)|= {n-1\choose m}{kn-1\choose m-1},
\end{equation}
\begin{equation}\label{c2}
|\mathcal{S}^{UD}_n(k, m)|= {n-1\choose m-1}{kn-1\choose m-1}.
\end{equation}
\end{lemma}
\pf Each $P\in \mathcal{S}^{UU}_n(k, m)$ can be  uniquely written as
$$UU^{x_1}D^{y_1}U^{x_2}D^{y_2}\ldots U^{x_m}D^{y_m}U^{x_{m+1}}$$
such that
$$
\left\{
\begin{array}{l}
 1+x_1+x_2+\ldots +x_{m+1}=n \\
 y_1+y_2+\ldots +y_{m}= kn,
\end{array}\right.
$$
where   $x_i, y_i\geq 1$  for $1\leq i\leq m$ and $x_{m+1}\geq 0$.
The solutions of $x_i$'s is equal to ${n-1\choose m}$ and the
solutions of $y_i$'s is equal to ${kn-1\choose m-1}$. Thus, Formula
(\ref{c1}) is proved.

Each $P\in \mathcal{S}^{UU}_n(k, m)$ can be uniquely written  as
$$U D^{y_1}U^{x_1}D^{y_2}\ldots U^{x_{m-1}}D^{y_m}U^{x_{m}}$$
such that
$$
\left\{
\begin{array}{l }
 1+x_1+x_2+\ldots x_{m}=n\\
 y_1+y_2+\ldots +y_{m}= kn,
\end{array}\right.
$$
where   $x_i,   y_i\geq 1$ for $1\leq i\leq m-1$, $y_m\geq 1$ and
$x_{m}\geq 0$. The solutions of $x_i$'s is equal to ${n-1\choose
m-1}$ and the solutions of $y_i$'s is equal to ${kn-1\choose m-1}$.
This leads to Formula (\ref{c2}).  This completes the proof. \qed

 From Formulae (\ref{c1}) and (\ref{c2}), we deduce that the
number of $k$-ary paths  with $n$ up steps and  $m$ peaks is equal
to
$$
{1\over m}({n-1\choose m-1}{kn-1\choose m-2}+{n-1\choose
m-1}{kn-1\choose m-1})={1\over n}{n\choose m}{kn\choose m-1}.
$$
Note that when $k=1$, we are led to the Narayana numbers \cite{De,
Na}.

 \noindent{\bf Acknowledgments.}   The   author was supported
by the National Natural Science Foundation of China (No.10901141).


 \newpage

\end{document}